\documentclass[10 pt,a4paper,twoside,reqno]{amsart} 
\usepackage{amsfonts,amssymb,amscd,amsmath,enumerate,verbatim,calc} 
\usepackage[dvipsnames,svgnames,x11names,hyperref]{xcolor}
\usepackage[a-1b]{pdfx}
\hypersetup{colorlinks,breaklinks,
             urlcolor=blue,
             linkcolor=blue}
 
\begin{document}
	
 \title[Role of History]{The Role of History in Learning and Teaching Mathematics: A Personal Perspective}

\thanks{This is a slightly revised version of the Presidential address (General) delivered 
at the 84th Annual Conference of the  Indian Mathematical Society held at Shri Mata Vaishno Devi University, Jammu, 
November 27-30,~2018. The author gratefully acknowledges comments and suggestions 
received from several colleagues, especially U. K. Anandavardhanan and Amartya K. Dutta.} 

\author{Sudhir R. Ghorpade}
	\address{Department of Mathematics, \newline \indent
		Indian Institute of Technology Bombay,\newline \indent
		Powai, Mumbai 400076, India.}
	\email{\href{mailto:srg@math.iitb.ac.in}{srg@math.iitb.ac.in}}
\urladdr{\href{http://www.math.iitb.ac.in/$\sim$srg/}{http://www.math.iitb.ac.in/$\sim$srg/}}

\maketitle

At the outset, I would like to express my gratitude to the Indian Mathematical Society for electing me as its President for the year 2018-19. I feel humbled by this honour and I consider it a privilege to address this august gathering of fellow mathematicians and invited guests. The Indian Mathematical Society, or in short, IMS, has a distinguished history of well over a century and many stalwarts have held the high office of the President of the IMS. Just to name a few, I might mention R. Vaidyanathaswamy, T. Vijayaraghavan, Hansraj Gupta, R. P. Bambah, R. B. Bapat, R. Sridharan, and S. G. Dani. This conference is rather unique because it is perhaps for the first time that it is scheduled  in the month of November and also the first time it is held in the state of Jammu and Kashmir. 

The topic I have chosen for this general talk is the role of history in learning and teaching mathematics. It would be beyond my competence to give here a comprehensive account of the history of mathematics or to give a scholarly narration befitting a professional historian. Instead, I will attempt to provide a personal perspective by relating some incidents in my own journey as a student and teacher of mathematics, especially those related with the history and development of the magnificent edifice that we call mathematics. Along the way, I will also mention some of my favourite books and articles. It may be hoped that this could be of some value and interest to younger students and researchers. I am reminded here of the following lines from a highly readable,  instructive, and MAA Chauvenet prize winning article \cite{AbH} 
by my teacher  Prof. Shreeram Abhyankar: 
\begin{quote}
Now it is true that the history of mathematics should primarily consist of an account of the achievements of great men; but it may be of some interest to also see their impact on an average student like me. With this in mind I may be allowed some personal reminiscences.
\end{quote}
I think Abhyankar was exceptionally humble to call himself \emph{an average student} here. But if I express a similar sentiment, then it would not be a matter of humility, but simply being truthful. Nonetheless, I hope you will permit some indulgence at this age and stage. I 
would also like 
to use the opportunity to pay a tribute to some of my teachers who 
played an important role in shaping my~career. 

Coming back to the theme of this address, the importance of history can 
hardly be overstated. In the preface to his 
book on differential equations, Simmons  \cite{Sim} mentions an old Armenian saying ``He who lacks a sense of the past is condemned to live in the narrow darkness of his own generation", and goes on to write:  
\begin{quote}
Mathematics without history is mathematics  stripped of its greatness:  for,  like the  other arts--and  mathematics  is one of  the supreme arts of civilization--it derives its grandeur from the fact of being a human creation.
\end{quote}
At any rate, this book of Simmons, and especially the historical notes in it, which I read during the second year of my undergraduate study had a considerable influence on me then. Moreover, it may have sown the seeds of eventually opting to make a career in mathematics. I was also blessed to have had very good college teachers like Prof. (Mrs.) S. B. Kulkarni and  Prof. V. R. Kulkarni. In the third year of my B.Sc., they encouraged me  to become in-charge of college's Mathematics Association and edit a wallpaper. This wallpaper was to contain snippets from the history of mathematics, problems, quotable quotes, etc. 
I remember writing several pieces under different names due to a lack of contributions from fellow~students. An old book of Cajori \cite{C} that I found in the college library was of great help to me then. This book has on its title page the following nice quote of J.\,W.\,L. Glaisher:
\begin{quote}
I am sure that no subject loses more than mathematics by any attempt to dissociate it from its history.
\end{quote}

For my Master's degree in Mathematics, I joined IIT Bombay in 1982 and this opened the doors to its well stacked library. It was quite exciting to be able to walk among books and browse through anything that seemed interesting. It was also thrilling to read the following lines by Sarvapalli Radhakrishnan that are engraved on a stone plaque in the foyer of the Convocation Hall of the Institute. 
\begin{quote}
Knowledge is not something to be packed away in some corner of our brain, but what enters into our being, colours the emotion, haunts our soul, and is as close to us as life itself. It is the overmastering power which through the intellect moulds the whole personality, trains the emotion and disciplines the will.
\end{quote}
At IIT Bombay, I  had the good fortune to be taught by some wonderful teachers such as Professors M. V. Deshpande, K. D. Joshi, B. V. Limaye, and D. V. Pai. In the summer vacation between the two years of my M.Sc., I had an opportunity to attend a summer school on Analysis and Probability conducted by the Indian Statistical Institute at the University of Mysore. Prior to going there, my friend Subathra Ramanathan gifted me with a copy of E. T. Bell's \emph{Men of Mathematics}. 
This book has 29 chapters, each giving an account of the life and work of a great mathematician in an absorbing and lucid manner. While I was in Mysore for a month (around the same time when India won the Cricket World Cup for the first time!), there would be a stern knock on the door very early in the morning to serve coffee, and after that one had to wait a considerable time for breakfast. During this period I would read one chapter each of Bell's book everyday,  and could thus  finish reading the book by the end of the month long school. The Mysore summer school was also 
memorable for me because of excellent teachers such as Professors B. V. Rao, S.  C. Bagchi, and V. S. Sunder, and brilliant fellow students like C. S. Rajan and V. Balaji. I was to learn later that although the historical account in Bell's book is interesting and readable, it wasn't always the most accurate. However, the book does give a broad perspective of the subject and some of its heros in a manner that can captivate an impressionable mind, and can still be recommended to young readers. But as indicated above, Bell's treatment, especially of ancient Indian mathematics and mathematicians, could be tempered with other sources. See, for example, Amartya Dutta's article \cite{DBh} on and in  \emph{Bh$\bar{\textit{a}}$van$\bar{\textit{a}}$}, and the references therein, as well as the article \cite{Sei} of Seidenberg (who is related to me by mathematical geneaology\footnote{One can trace mathematical geneaology at: 
\href{https://www.genealogy.math.ndsu.nodak.edu/}{https://www.genealogy.math.ndsu.nodak.edu/}}), 
which is also a good reference for  real ``Vedic Mathematics". Furthermore, today's reader should complement Bell \cite{Bel} with a book such as \cite{WM} on women of mathematics, and read something (for instance,  \cite{JA}), about the first female Fields medalist Maryam Mirzakhani. 

It was also in the year 1983 that I first met Prof. Shreeram Abhyankar at Bombay University (rather accidentally, as I have mentioned in \cite{Gh1}) and later went on to do Ph.D. under his supervision at Purdue University. In the course of my long association with Prof. Abhyankar, I learned a great deal not only about mathematics and many mathematicians he knew personally, but also about his view of the historical development of mathematics, 
especially, algebra and algebraic geometry. Abhyankar spoke highly of 
Felix Klein's \emph{Entwicklung}, and while I couldn't read the German original, I did read parts of the English translation \cite{FK} and could see what a valuable resource it was. Somewhere around this time, probably much earlier, I came across G. H. Hardy's \emph{Apology} \cite{H}. What a remarkable book this is! It gives a brilliant insight into the mind of a mathematician. Even if some of what Hardy writes is controversial and perhaps untrue, I would heartily recommend this book to anyone interested in mathematics, or for that matter, any layman curious about mathematics and  mathematicians. I certainly find myself going back to this book from time to time, even if just to admire Hardy's turn of phrase.  
For a more in-depth look into the working of a mathematician's mind, Hadamard's essay \cite{Had} (which I have only skimmed through) could be consulted. Perhaps a more riveting account is given by Poincar\'e in \cite{P}
where, among other things, he narrates how the crucial idea about automorphic functions came to him just as he was about to put his foot on a bus. The book also  has a scholarly preface by Bertrand Russell, which makes for an interesting~read.  In his 
 ICM address \cite{W1}, Andr\'e Weil highlights the importance of history for creative and would-be creative scientists, and includes (an English translation of) the following words of Leibniz:
\begin{quote}
Its use is not just that History may give everyone his due and that others may look forward to similar praise, but also that the art of discovery be promoted and its methods known through illustrious examples. 
\end{quote}

Biographies of great mathematicians can also teach us a lot and help gain a perspective on their work. There are several that are available, but I will mention just one that I remember reading with pleasure. Namely, Constance Reid's biography of Hilbert \cite{CH}, which Freeman Dyson called ``a poem in praise of mathematics."  Another inspiring book of a more recent vintage is Sreenivasan's compilation \cite{100} of accounts by about 100 distinguished scientists, including several mathematicians, that outline their motivations for taking up science and a brief description of their contributions to the subject. I would add to this collection a fascinating and instructive interview \cite{Ser} with one of the greatest mathematicians of our time, J.-P. Serre, whom I have had the pleasure of meeting at Purdue, Paris and Luminy. 

While biographical details and anecdotes about the life and times of mathematicians can be interesting, a mathematician finds true expression in his or her mathematical work. Thus for a serious student of mathematics, there is often no better alternative to learn the subject than to go through original sources. Isaac Newton famously spoke of being able to see farther because he was standing on the shoulders of giants. 
Abel puts this point across more directly when he says: 
\begin{quote}
... if one wants to make progress in mathematics one should study the masters not the pupils. 
\end{quote}
This particular quote appears at the top of the translators' preface to Hecke's \emph{Vorlesungen} \cite{Hec}, which is a masterly treatment of the theory of algebraic numbers. Speaking of algebraic number theory, Hilbert's \emph{Zahlbericht} \cite{Hil} is now available in English and it is still an excellent introduction to the subject. As for elementary number theory, it is still profitable to read Gauss's \emph{Disquisitiones} \cite{Gau}, which was first published in Latin in 1801 and is available in English. Many good libraries contain collected works of eminent mathematicians, and they make it easier for us to understand the evolution of ideas of these mathematicians and the paths taken in their mathematical journeys. Today, thanks to the Internet, it has become simpler to access many original sources and it may suffice to 
cite the Digital Mathematics Library (\href{https://www.math.uni-bielefeld.de/$\sim$rehmann/DML/}{https://www.math.uni-bielefeld.de/$\sim$rehmann/DML/}) and the 
AMS Digital Mathematics Registry (\href{https://mathscinet.ams.org/dmr/}{https://mathscinet.ams.org/dmr/}). 

Let me turn to the role of history in teaching mathematics. Here too, rather than pontificating about generalities, I would like to draw from personal experiences. By the end of 1989, I returned to IIT Bombay to join as a faculty member, and taught for the first time a  course on algebra for M.Sc. students in Spring~1990. (Incidentally, another thing I did in 1990 was to become a life member of the Indian Mathematical Society, and to my surprise, I was invited to give a talk at the annual conference of the IMS at Surat later that year.) Partly because of the course I was teaching and partly due to my own interest, I consulted the book of van der Waerden on history of algebra \cite{vdW}. It was then that I understood why abelian groups were named after Abel. This book also helped me to  understand better the evolution of theory of equations and the advents made by 
Scipione del Ferro, Tartaglia, Cardano, Ferrari, and how they culminated in the monumental work of  \'Evariste Galois (1811-1832) which brought to fore the basic notion of a group. Later I  could use this (I think with good effect) in lectures for younger audiences that explained the development of algebra. Over the years, my hobby of picking up bits and pieces of history of mathematics has, I believe, stood me in good stead, and perhaps enhanced the quality of my teaching. Students 
are often amused by snippets of history such as 
for instance, the fact that the basic notion of an ideal in a ring came about in Kummer's attempt to prove Fermat's Last Theorem, or that many of the fundamental results of commutative algebra and algebraic geometry (such as Hilbert's basis theorem,  Nullstellensatz, and syzygy theorem) were, in fact, lemmas in Hilbert's breakthrough work on invariant theory. Even in a very basic course such as Calculus, students are sometimes surprised to learn that the development of calculus did not happen in the order it is given in their textbooks. The notion of integration came much before the notion of differentiation was discovered! Realizing this, one can perhaps have a much better appreciation of the fundamental theorem of calculus, which relates seemingly disparate geometric problems of finding area under a curve in the plane and determining the (slope of) tangents to plane curves. 
Once again, thanks to the Internet, it is no longer necessary for an interested student to leaf through dusty volumes in the library, but simply access the requisite information on a digital media of his or her choice in order to get a fairly good sense of the history of a topic in mathematics or a mathematician. For them, I would recommend using the MacTutor History of Mathematics archive (\href{http://www-history.mcs.st-and.ac.uk/}{http://www-history.mcs.st-and.ac.uk/}). Many a time, the Wikipedia is also a good source of information. The historical notes appearing at the end of Bourbaki's \emph{Elements} 
have long been recognized as 
an authoritative source. A compilation of most of these historical notes is available in book form \cite{B}. 
A very nice collection of articles about American mathematics can be found in the first 3 volumes of the AMS series on History of Mathematics \cite{D}. 
I should also 
mention a few good sources about mathematics in India and 
Indian mathematicians. 
Here. I warmly recommend  Volumes 3  and 5 of the \emph{Culture and History of Mathematics} series of HBA  \cite{ESS, Se}, and the scholarly articles of Raghavan Narasimhan \cite{RN} and M. S. Raghunathan \cite{MSR}. 

Let me end with a somewhat philosophical remark. Why do we do mathematics? And why should we learn something about its history? The first question is relatively easy to answer: we do mathematics because it brings us joy (not to mention a job and a salary ...). I have dwelt on this question in greater details elsewhere\footnote{\emph{Empowering Times}, Jan 2018. \href{http://empoweredindia.com/Archives\_ET/ET201801.html}{http://empoweredindia.com/Archives\_ET/ET201801.html}}, quoting  Hardy \cite{H} rather extensively, and so I will not discuss it further, except to say that we mathematicians yearn for those rare moments of exhilaration when suddenly we understand some seemingly complex topic very clearly, with everything falling into place, or, if we are luckier still, then we discover something new and nice. In a way, this brings us closer to divinity, albeit for a brief fleeting moment. Needless to say, for first rate mathematicians, such moments of exhilaration would be more frequent and of a much higher order.  By dwelling upon the history of mathematics and the achievements of distinguished mathematicians, we 
not only see glimpses of greatness, but an easier way 
to be a step closer to divinity, even if vicariously. 
We also feel a greater connect with the mathematical milieu. At any rate, whatever be your reason, I wish you many pleasant 
moments with mathematics and its history, and an enjoyable conference!

\end{document}